\documentclass[reqno]{amsart}
\usepackage{graphicx}
\usepackage{amsmath}
\usepackage{amsthm}
\usepackage{mathtools}
\usepackage{amsfonts}
\usepackage{amssymb}
\usepackage{faktor}
\usepackage{wrapfig}
\usepackage{epstopdf}
\usepackage{caption}
\usepackage{url}
\usepackage{mathtools} 
\usepackage{float}
\usepackage{caption}
\usepackage{centernot}
\usepackage{subcaption}
\usepackage{accents}

\title{
Reducible surgery in lens spaces and seiferters}

\author{Fyodor Gainullin}

\newtheorem{theorem}{Theorem}

\newtheorem{prop}[theorem]{Proposition}

\newtheorem{qn}[theorem]{Question}

\newtheorem{conj}[theorem]{Conjecture}

\begin{document}

\begin{abstract}
The Cabling Conjecture states that surgery on hyperbolic knots in $S^3$ never produces reducible manifolds. In contrast, there do exist hyperbolic knots in some lens spaces with non-prime surgeries. Baker constructed a family of such hyperbolic knots in \cite{bakerCabling} and posed a conjecture that his examples encompass all hyperbolic knots in lens spaces with non-prime surgeries. Using the idea of seiferters (see \cite{networking1}) we construct a counterexample to this conjecture. In the process of construction, we also derive an obstruction for a small Seifert fibred space to be obtainable by a surgery with a seiferter.
\end{abstract}

\maketitle

\section{Introduction}
\label{sec:intro}

Given a hyperbolic knot, if we look at all manifolds we can produce by surgery on this knot, it will turn out that only finitely many of them are not hyperbolic. In fact, merely saying `finitely many' doesn't do justice to the result -- there are at most 10 of these so-called \slshape exceptional \upshape surgeries on every hyperbolic knot \cite{lackenbyMeyerhoff}.

This naturally leads to the feeling that understanding the reason for these exceptions will help us understand the structure of three-dimensional manifolds and knots in them. Indeed, understanding exceptional surgeries is a major effort in contemporary 3-manifold topology.

Results of exceptional surgery in $S^3$ are either non-prime, Seifert fibred or toroidal (see \cite{gordon} and references therein). The Cabling Conjecture \cite{gonzalez-acunaShort} states that exceptional surgery, in fact, never produces non-prime manifolds. There has been a lot of progress in this direction: we know that the slope has to be integral \cite{gordonLueckeReducible}, the result of reducible surgery on any knot always contains a lens space summand \cite{gordonLueckeKnotComplement} and the cabling conjecture is true when the surgery results in a connected sum of lens spaces \cite{greeneCabling}. It has also been proven for large classes of knots (see \cite{menascoThistlethwaite, eudave-munozCabling, hayashiShimokawa, wuArborescent, hoffmanThesis, scharlemannReducible}).

In contrast to the reducible case, there are numerous examples of exceptional surgeries yielding Seifert fibred manifolds \cite{berge, bleilerHodgson, boyerZhang, brittenhamWu, dean, eudave-munoz, teragaito, mattmanMiyazakiMotegi, networking1}. Non-toroidal Seifert fibred spaces that can be obtained by surgery on knots are subdivided into two families: lens spaces and small Seifert fibred spaces with three exceptional fibres.  We seem to know considerably less about the second case: the integrality of slope is only conjectured \cite{gordon2} and there is no proposed list of knots that admit such surgeries. However, there is a conjectural reason for such surgeries \cite{networking1}: most knots with Seifert fibred surgeries are known to be related to torus knots by twisting along seiferters (see below for the definition).

A first slight generalisation of $S^3$ is a lens space. One might wonder if a statement as strong as the Cabling conjecture holds for general lens spaces. This turns out to be false -- examples of hyperbolic knots with surgeries yielding connected sums of lens spaces have been found in \cite{boyerZhangSeminorms, eudave-munozWu, kang, bakerCabling}. Baker \cite{bakerCabling} noticed that all these examples can be explained by a single unifying construction based on rational tangle replacement in a certain family of 3-braids. On the basis of this observation he put forward the following conjecture \footnote{Baker's conjecture is in fact stronger -- it also lists the exact situations when such surgeries occur. This statement, however, is enough for our purposes.}.

\begin{conj}[Baker]
Assume a knot $K$ in a lens space admits a surgery to a non-prime 3-manifold $Y$. If $K$ is hyperbolic, then $Y = L(r, 1) \# L(s, 1)$. Otherwise either $K$ is a torus knot, a Klein bottle knot, or a cabled knot and the surgery is along the boundary slope of an essential annulus in the exterior of $K$, or $K$ is contained in a ball.
\label{baker_conj}
\end{conj}

In fact, Baker proves the part of his conjecture concerning non-hyperbolic knots. In this paper we provide an example that does not fit into this construction\footnote{We have learned from personal correspondence with Ken Baker that he has also produced examples contradicting his conjecture. His construction appears to be very different from ours.}.

\begin{theorem}
There is a hyperbolic null-homologous knot $K' \subset L(15,4)$ of genus $1$ that gives $L(5,3)\# L(3,2)$ by surgery.
\label{main1}
\end{theorem}

Our construction is based on the idea of seiferters which were defined in \cite{networking1}. We obtain a knot in a lens space with a Seifert fibred surgery by first doing a surgery on the seiferter. To know which lens space we get as a result we first need to consider how the preimages of the other fibres wrap around our chosen seiferter. This turns out to depend on the square of the linking number of the sieferter with the knot, which leads to Theorem \ref{main2} below. The idea of using seiferters for a counterexample led us to consider many examples of knots with Seifert fibred surgeries and somewhat surprisingly most of them did not give a counterexample we were after. During this search we made some observations about currently known examples of Seifert fibred surgeries with seiferters. One of such observations is summarised in Proposition \ref{berge} below.

\begin{figure}
\includegraphics[scale=0.32, clip = true, trim = 35 495 35 40]{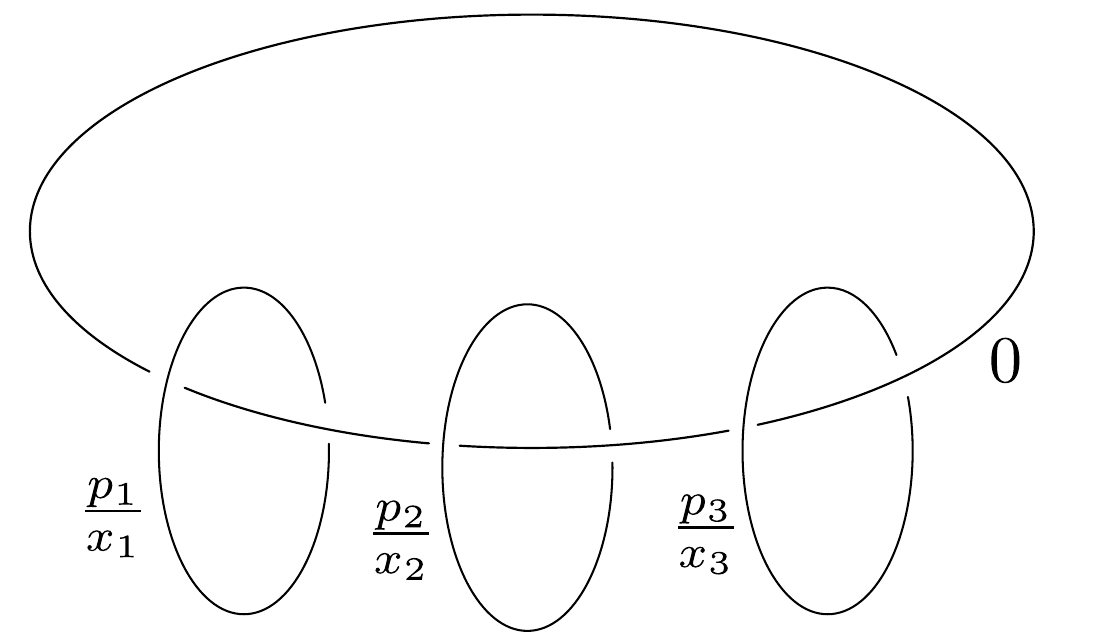}
\caption{Surgery description for $S^2((p_1,x_1), (p_2,x_2), (p_3, x_3))$.}
\label{SFSdefn}
\end{figure}

In Theorem \ref{main2}, $S^2((p_1,x_1), (p_2,x_2), (p_3, x_3))$ denotes the space with the surgery presentation as in Figure \ref{SFSdefn}. It is a Seifert fibred space over a sphere with three exceptional fibres (with multiplicities $p_1, p_2$ and $p_3$, which are positive integers).

\begin{theorem}
Suppose $S^2((p_1,x_1), (p_2,x_2), (p_3, x_3))$ is obtained by a surgery on a knot in $S^3$ with a seiferter. Let $H = p_1p_2x_3+p_1p_3x_2 + p_2p_3x_1$. For $i = 1, 2, 3$ let $q_i$ be a multiplicative inverse of $x_i$ modulo $p_i$. If $H \neq 0$ then one of
\begin{equation}
\pm \frac{q_1H-p_2p_3}{p_1}, \ \pm \frac{q_2H-p_3p_1}{p_2}, \ \pm \frac{q_3H-p_1p_2}{p_3}, \ \pm p_1p_2p_3
\label{obstructRatHom}
\end{equation}
is a quadratic residue modulo $H$ (for some choice of sign).

If $H = 0$, then
\begin{equation}
p_1p_2p_3
\end{equation}
is a square.
\label{main2}
\end{theorem}

We remark that up to sign, $H$ is the order of the first homology of \linebreak $S^2((p_1,x_1), (p_2,x_2), (p_3, x_3))$ (we interpret infinite homology groups as having order $0$, so $H = 0$ means infinite first homology group).

Proposition \ref{prop:seiferter_restrictions} is a more general surgery restriction which directly implies Theorem \ref{main2}. Thus in Section 3 we prove Proposition \ref{prop:seiferter_restrictions}, which implies Theorem \ref{main2}.

Note that if $(p_i, p_j) = 1$ for all $i\neq j$, then condition \eqref{obstructRatHom} is equivalent to $\pm p_1p_2p_3$ being a quadratic residue.

Even though Theorem \ref{main2} is vacuous when $H$ is a prime congruent to $3$ modulo $4$, there do exist many examples when it gives an obstruction for a space to be obtained by Seifert surgery with a seiferter. Unfortunately, Theorem \ref{main2} also does not give any information about the examples of Teragaito from \cite{teragaito}. The next proposition gives an example of such an infinite family.

\begin{prop}
Let $p \equiv 3 \ (\mbox{mod } 17)$. Then $S^2((p,-17),(2p-1,17),(2p+1,17))$ cannot be obtained by a Seifert surgery on a knot in $S^3$ with a seiferter.
\label{obstructExamples}
\end{prop}

We do not know if the spaces from the proposition above can be obtained by an integral surgery on a knot in $S^3$.

Finally, we notice that often there exists a surgery on a seiferter that makes the original knot with a Seifert fibred surgery non-hyperbolic in the resulting lens space. This holds, in particular, for all Berge knots so we have (we consider a core of one of the solid tori to be a torus knot for the statement below)

\begin{prop}
Every Berge knot can be obtained from a torus knot or a cable of a torus knot by repeatedly taking band sums with a cable of some unknot. For each Berge knot the unknot and its cable are fixed.
\label{berge}
\end{prop}

\section*{Acknowledgements}

I am extremely grateful to Ken Baker for many useful discussions and comments he made about the paper. He also found the filling one needs to perform on the seiferter of Figure \ref{pretzel4} to turn the knot into a torus knot. I would also like to thank my supervisor Dorothy Buck for her continued support over the course of my Ph.D. studies.

\section{Overview of the construction and some obvious examples that do not work}
\label{sec:overview}

For the current paper we define the lens space $L(p,q)$ to be the result of $-p/q$-surgery on the unknot. Denote by $K(m)$ the result of $m$-surgery on a knot $K$.

Suppose $K \subset S^3$ is a knot with a Seifert fibred surgery. Assume additionally that there is an unknotted simple closed curve $c$ disjoint from $K$ which becomes one of the Seifert fibres after the surgery. Such a curve is called a \slshape seiferter\upshape. If $K$ is not a torus knot and such $c$ exists, then the surgery is integral and conjecturally every integral Seifert fibred surgery possesses a seiferter \cite{miyazakiMotegiIII}.

If $c$ is a seiferter with Seifert invariant $(p,x)$, then in the Seifert fibred solid torus neighbourhood of $c$ other fibres are $(p,q)$-curves where $q$ is some multiplicative inverse of $x$ modulo $p$.

In the solid torus complementary to a neighbourhood of $c$ the fibres on the boundary are $(q,p)$-curves. So if we perform a surgery on $c$ with the slope given by the fibres around it, we get the lens space $L(q,p)$.

Denote by $K'$ the image of the knot $K$ in this lens space. Suppose we now perform the surgery on $K'$ with the slope induced by the original surgery slope on $K$. The result is the same as first performing the original Seifert fibred surgery on $K$ and then doing the reducible surgery on the fibre that $c$ becomes.

Drilling out a fibre $h$ with a Seifert invariant $(p_n, x_n)$ from the Seifert fibred space $S^2((p_1,x_1), (p_2,x_2), \ldots (p_n, x_n))$ and regluing a solid torus with the slope given by ordinary fibres on the boundary of the neighbourhood of $h$ produces $\#_{i = 1}^{n-1} L(p_i,x_i)$ \cite{heil}.

So if the surgery on $K$ produced $S^2((p_1,x_1),(p_2,x_2),(p_3,x_3))$ and $c$ became the exceptional fibre of index $p_3$, then this construction gives a knot in a lens space $L(q_3,p_3)$ with a surgery yielding $L(p_1,x_1)\#L(p_2,x_2)$. Here $q_3$ is some inverse of $x_3$ modulo $p_3$, the exact value of which we find in the next section.

Now in order to disprove Conjecture \ref{baker_conj} it is enough to find a knot $K \subset S^3$ with a surgery yielding $S^2((p_1,x_1),(p_2,x_2),(p_3,x_3))$ such that the following conditions hold:
\begin{itemize}
\item there is a seiferter $c$ for this surgery that becomes the exceptional fibre of index $p_3$;
\item $x_1 \not \equiv \pm 1$ modulo $p_1$;
\item after the appropriate surgery on $c$ the knot $K'$ (the image of $K$) is hyperbolic in the resulting lens space.
\end{itemize}

For a Seifert invariant $(p, x)$ of a Seifert fibre, the first number in the pair, $p$, is known as its multiplicity. However, it seems there is a lack of a widespread term for the second number, $x$. In fact, it is only defined modulo the multiplicity. In this paper, we call this congruence class of $x$ the \slshape torque \upshape of the Seifert fibre. We will slightly abuse this terminology by calling any representative of the congruence class of $x$ the torque of the Seifert fibre.

As is clear from our strategy, we are searching for examples of surgeries on knots in $S^3$ that produce Seifert fibred spaces over the sphere for which some exceptional fibres have torque not equal to $\pm 1$. \footnote{In a previous version of this paper we asked a question whether one of the torques always has to be $\pm 1$ modulo the corresponding fibre multiplicity, as all examples we encountered did have this property. It was pointed out by John Berge, however, that examples that do not satisfy this exist.}

The obvious big collections of small Seifert fibred surgeries with seiferters to consider in our search are those of Eudave-Mu\~noz \cite{eudave-munoz} and Berge \cite{berge}. The seiferters for these are exhibited in \cite{networking4} and \cite{networking2} respectively.

Now we need to take an example of a surgery with a seiferter such that at least one of the exceptional fibres which is not the image of our seiferter has torque not equal to $\pm 1$. In the examples of \cite{eudave-munoz} it turned out that if at least two torques were $\pm 1$ then both such fibres were not in the image of our seiferter. Our initial attempt was to take one of the infinite families in \cite{eudave-munoz} in which only one of the torques is $\pm 1$ and turn it into an infinite family of hyperbolic knots in lens spaces that give connected sums of lens spaces by surgery. However, this attempt failed -- the knots that we obtained turned out to be torus knots or cables of torus knots.

This is somewhat surprising due to two facts. Firstly, the knots we are dealing with are hyperbolic in $S^3$ and seiferters become exceptional fibres. Secondly, due to \cite[Corollary 3.14]{networking1} if a hyperbolic knot has a seiferter that becomes an exceptional fibre, then the link formed by the knot and the seiferter is hyperbolic. Thus `generically' one expects surgeries on such seiferters to give hyperbolic knots.

A similar phenomenon occurs for Berge knots. We illustrate this for Berge knots of types IX-X by a sequence of figures. 

\begin{figure}
\includegraphics[scale=0.5, clip = true, trim = 5 370 95 170]{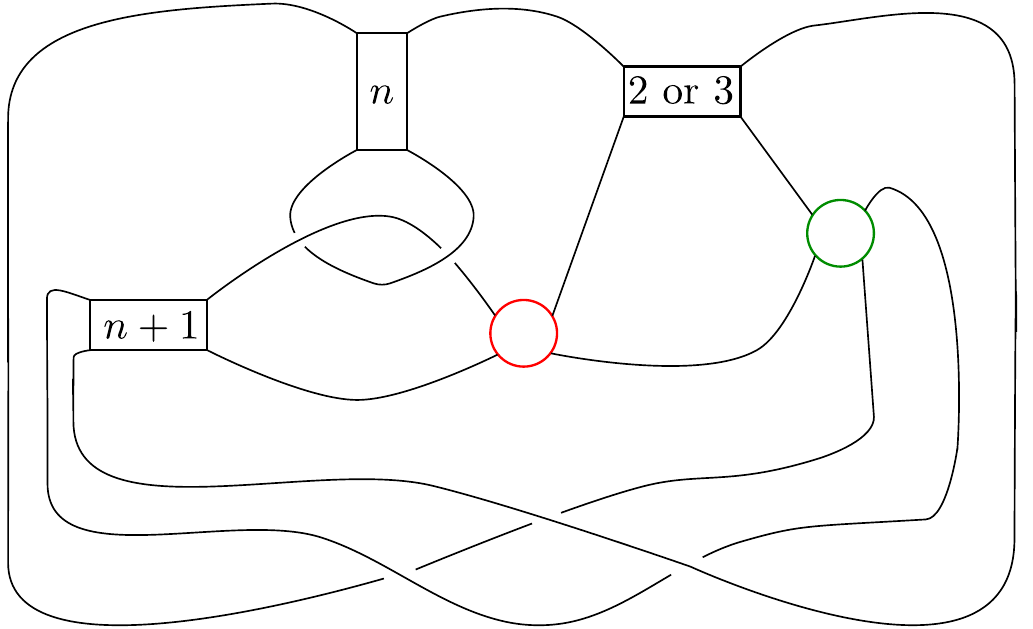}
\caption{Tangle description for the links formed by Berge knots of type IX-X and seiferters. The knot component is in red and the seiferter component is in green.}
\label{bergeSporadicAB1}
\end{figure}

Consider the tangle in Figure \ref{bergeSporadicAB1}. This is the same tangle as \cite[Figure 41]{bakerBerge2} with one more tangle (in green) removed (we use conventions of \cite{networking4} for the diagrams of tangles). In other words, there is a rational filling of the green tangle that makes the double branched cover of the resulting tangle the exterior of a Berge knot of type IX or X. Filling the red component with $\infty$-tangle (which corresponds to the trivial surgery on the corresponding Berge knots) makes this a trivial tangle, thus the ball bounded by the green sphere lifts to a Heegaard solid torus of $S^3$, i.e. its core is unknotted in $S^3$. The lens space surgery corresponds to filling the red component with the $0$-tangle. The tangle then becomes a sum of two rational tangles so we can see that the unknot corresponding to the green component in the Figure \ref{bergeSporadicAB1} does become a Seifert fibre. Thus it is a seiferter.

\begin{figure}
\includegraphics[scale=0.5, clip = true, trim = 30 405 160 170]{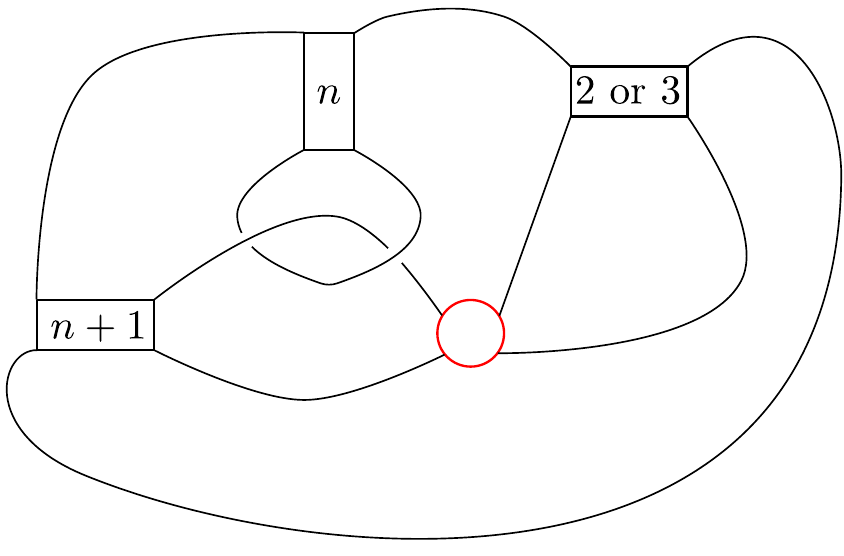}
\caption{Tangle description of the images of Berge knots of types IX-X after a surgery on a seiferter.}
\label{bergeSporadicAB3}
\end{figure}

We now want to show that there is a surgery on a seiferter which makes the original Berge knot a cable of a torus knot. For this end, fill the green component with the $\infty$-tangle. The tangle then becomes as in Figure \ref{bergeSporadicAB3}.

\begin{figure}
\includegraphics[scale=0.5, clip = true, trim = 0 630 150 5]{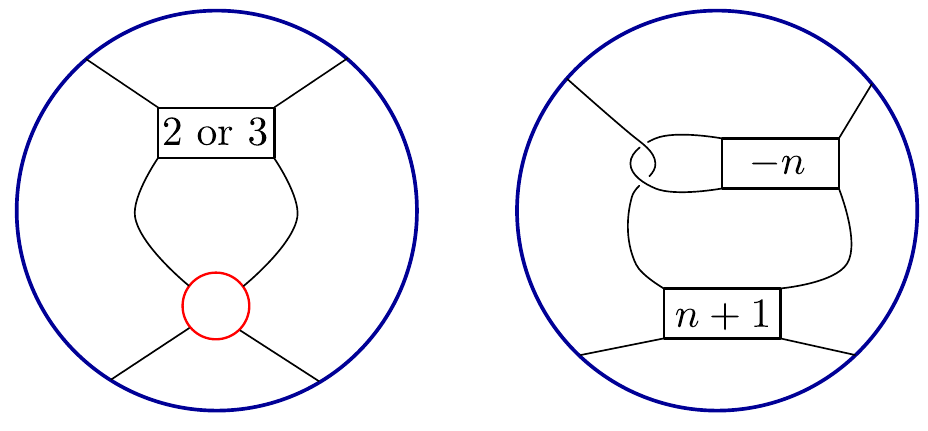}
\caption{The tangle of Figure \ref{bergeSporadicAB3} decomposes into a union of two tangles.}
\label{bergeSporadicAB4}
\end{figure}

This decomposes into a union of tangles as in Figure \ref{bergeSporadicAB4}. The first of these tangles (from the left) is a cable space and the second is a torus knot exterior (in a lens space), thus the knot is a cable of a torus knot.

A similar argument works for all other Berge knots apart from those of types VII and VIII and for knots of Eudave-Mu\~noz from \cite{eudave-munoz}. We will illustrate this for other Berge knots later.

Unfortunately, this stopped us from producing infinite family of counterexamples. In the end we did succeed by finding a counterexample based on the Seifert fibred surgery from \cite{mattmanMiyazakiMotegi}.

Now to know which lens space we want to end up in we need to find how preimages of the ordinary fibres wind around the seiferter in $S^3$. This is the purpose of the next section.

\section{Fibred neighbourhood of a seiferter in $S^3$ and the proof of Theorem \ref{main2}}
\label{sec:fibred_nbrhd}

The aim of this section is twofold. As indicated in the previous section, we want to understand what the preimage of the fibred solid torus neighbourhood of a seiferter looks like in $S^3$. More precisely, we want to know how many times the ordinary fibres around the seiferter wind in each direction. Finding an expression that describes this winding in terms of seifert invariants and the linking number of the seiferter with the knot is the first aim of this section. As a result of this expression we also prove Theorem \ref{main2}, which is the second aim of this section.

\subsection{Exceptional fibre}

Suppose that $m$-surgery on a knot $K$ in $S^3$ produces $S^2((p_1,x_1),(p_2,x_2),(p_3,x_3))$ and that this surgery has a seiferter which becomes the exceptional fibre with the Seifert invariant $(p_1, x_1)$. The first homology group of $S^2((p_1,x_1),(p_2,x_2),(p_3,x_3))$ has the following abelian presentation:
$$
<x,y,z,h|p_1x+x_1h=0, p_2y+x_2h=0, p_3z+x_3h=0, x+y+z=0>.
$$

Thus its relation matrix is given by
$$
\begin{pmatrix} p_1 & 0 & 0 & x_1 \\ 0 & p_2 & 0 & x_2 \\ 0 & 0 & p_3 & x_3 \\ 1 & 1 & 1 & 0 \\  \end{pmatrix}.
$$

Modulus of the determinant of a relation matrix of an abelian group is its order (to avoid special cases we say that an abelian group `has order 0' to mean that it is infinite). Let $H = p_1p_2x_3+p_1p_3x_2 + p_2p_3x_1$ and $\delta = sign(H)$. Note that $H$ does not depend on the particular representation of a small Seifert fibred space. Then by calculating the determinant of the above matrix we conclude that the order of the first homology of $S^2((p_1,x_1),(p_2,x_2),(p_3,x_3))$ is $\delta H$. Suppose for now that $H \neq 0$.

We also assume that $K(m) = S^2((p_1,x_1),(p_2,x_2),(p_3,x_3))$, so if $\epsilon = sign(m)$, we have
$$
m = \epsilon \delta H.
$$

Suppose we do $t$ twists along the seiferter and the linking number of the seiferter with the knot is $l$. Recall that the ordinary fibres around the seiferter are $(p_1,q_1)$-curves, where $q_1$ is a multiplicative inverse of $x_1$ modulo $p_1$. We want to find the value of $q_1$. According to \cite[Section 5.1]{networking1} $m+tl^2$ surgery along the resulting knot gives $S^2((tq_1 + p_1,t\frac{q_1x_1-1}{p_1}+x_1),(p_2,x_2),(p_3,x_3))$. If we choose $t$ to have the same sign as $m$, then we have 
$$
m+tl^2 = \epsilon | (tq_1+p_1)p_2x_3+(tq_1+p_1)p_3x_2+p_2p_3(t\frac{q_1x_1-1}{p_1}+x_1)| = \epsilon | H + t\frac{q_1H-p_2p_3}{p_1}|.
$$

Note that since we can change $t$ arbitrarily, $l = 0$ if and only if $q_1H = p_2p_3$. Suppose $l \neq 0$.

Let $\gamma = sign(\frac{q_1H-p_2p_3}{p_1})$. Then for $|t|$ large enough, $sign(H + t\frac{q_1H-p_2p_3}{p_1}) = sign(t\frac{q_1H-p_2p_3}{p_1}) = \epsilon \gamma$.

Therefore
$$
m+tl^2 = \gamma(H + t\frac{q_1H-p_2p_3}{p_1}),
$$
so together with $m = \epsilon \delta H$ this gives
$$
(\epsilon \delta - \gamma)H = t(\gamma\frac{q_1H-p_2p_3}{p_1}-l^2).
$$

We remind that $\frac{q_1H-p_2p_3}{p_1}$ is an integer, so we must have that for all big enough $t$, $t|(\epsilon \delta - \gamma)H$. Clearly, this means that $\epsilon \delta = \gamma$ and
\begin{equation}
\label{excep}
l^2 = \epsilon \delta \frac{q_1H-p_2p_3}{p_1} \mbox{ that implies } q_1 = \frac{\epsilon \delta p_1l^2+p_2p_3}{H} \mbox{ if } H \neq 0.
\end{equation}

Note that this equation also holds in the case $l=0$.

If $H=0$ a similar but simpler argument implies
\begin{equation}
l^2 = \frac{p_2p_3}{p_1}.
\label{exch0}
\end{equation}

However, $0 = H = p_1p_2x_3+p_1p_3x_2 + p_2p_3x_1$ and $p_i$ and $x_i$ are coprime, hence $p_1|p_2p_3$ automatically so the new information we obtain is equivalent to $p_1p_2p_3$ being a square.

In particular, $l$ cannot be equal to $0$ when $H = 0$.

\subsection{Ordinary fibre}

Suppose now that the seiferter is an ordinary fibre such that other ordinary fibres in its solid torus neighbourhood are $(1,n)$-curves. Before twisting, the order of the first homology of the resulting manifold is still $\epsilon m = \delta H$. Let $H \neq 0$. As before, if the seifert invariant of the seiferter was $(1,0)$, after twisting $t$ times it changes to $(tn+1,-t)$.

This means that the resulting manifold changes to
$$
S^2((p_1,x_1),(p_2,x_2),(p_3,x_3),(tn+1,-t)).
$$

If we set $sign(t) = sign(m) = \epsilon$, order of its first homology is
\begin{multline*}
\epsilon(m+tl^2) = | p_1p_2p_3(-t)+p_1p_2(tn+1)x_3+p_1p_3(tn+1)x_1+p_2p_3(tn+1)x_1 | = \\
 = | H+t(nH-p_1p_2p_3)|.
\end{multline*}

As before, $l = 0$ if and only if $nH = p_1p_2p_3$. Suppose $l \neq 0$.

Let $\gamma = sign(nH-p_1p_2p_3)$. Then for large enough $|t|$, we have
$$
\epsilon(m+tl^2) = \epsilon \gamma (H+t(nH-p_1p_2p_3)), \mbox{ so } \epsilon \delta H + tl^2 = \gamma H + \gamma t(nH-p_1p_2p_3)
$$
and it follows similarly to the previous case that $\gamma = \epsilon \delta$ and
\begin{equation}
\label{ord}
l^2 = \epsilon \delta (nH-p_1p_2p_3) \mbox{ that implies } n = \frac{\epsilon \delta l^2 + p_1p_2p_3}{H} \mbox{ if } H \neq 0.
\end{equation}

This equation also holds when $l = 0$.

If $H = 0$ we get
\begin{equation}
l^2 = p_1p_2p_3
\label{ordh0}
\end{equation}

As in the previous case, when $H = 0$ we also get that $l \neq 0$.

\subsection{Proof of Theorem \ref{main2}}
Note that equations \eqref{excep}, \eqref{exch0}, \eqref{ord}  and \eqref{ordh0} give obstructions to fibres being preimages of seiferters. Since $q_1$ is well-defined modulo $p_1$ and changing $q_1$ by multiples of $p_1$ changes the overall expression by multiples of $H$, \eqref{excep} and \eqref{ord} are well-defined modulo $H$. Hence we obtain

\begin{prop}
Suppose $S^2((p_1, x_1),(p_2, x_2),(p_3, x_3))$ is obtained by a surgery on a knot $K \subset S^3$ with a seiferter $c$. Let $q_i$ be an inverse of $x_i$ modulo $p_i$ for $i = 1, 2, 3$. If $c$ becomes an exceptional fibre with Seifert invariant $(p_1,x_1)$ and $H \neq 0$ we have
\begin{itemize}
\item $\delta \frac{q_1H-p_2p_3}{p_1}$ is a quadratic residue modulo $H$ if the surgery slope is positive;
\item $\delta \frac{p_2p_3-q_1H}{p_1}$ is a quadratic residue modulo $H$ if the surgery slope is negative.
\end{itemize}
If $c$ becomes an ordinary fibre, we have
\begin{itemize}
\item $-\delta p_1p_2p_3$ is a quadratic residue modulo $H$ if the surgery slope is positive;
\item $\delta p_1p_2p_3$ is a quadratic residue modulo $H$ if the surgery slope is negative.
\end{itemize}
If $H = 0$, $p_1p_2p_3$ is a square.
\label{prop:seiferter_restrictions}
\end{prop}

This proposition is just a slightly expanded version of Theorem \ref{main2}, so the theorem follows as well.

We are now ready to describe our counterexample to Conjecture \ref{baker_conj}.

\section{Proof of Theorem \ref{main1}}
\label{sec:main_thm}

{
\renewcommand{\thetheorem}{\ref{main1}}
\begin{theorem}
There is a hyperbolic null-homologous knot $K' \subset L(15,4)$ of genus $1$ that gives $L(5,3)\# L(3,2)$ by surgery.
\end{theorem}
\addtocounter{theorem}{-1}
}

\begin{proof}

Let $K$ be the $(-3,3,5)$-pretzel knot depicted in Figure \ref{pretzel}. Then the following is proven in \cite{mattmanMiyazakiMotegi}.

\begin{theorem}
Surgery with slope $1$ on $K$ gives the small Seifert fibred space \linebreak $S^2((5, -2),(3, -1),(4, 3))$. There is a seiferter $c$ for this surgery that becomes the exceptional fibre of index $4$. Moreover, there is a genus $1$ Seifert surface for $K$ that does not intersect $c$.
\end{theorem}

\begin{figure}
\includegraphics[scale=0.3, clip = true, trim = 85 160 105 150]{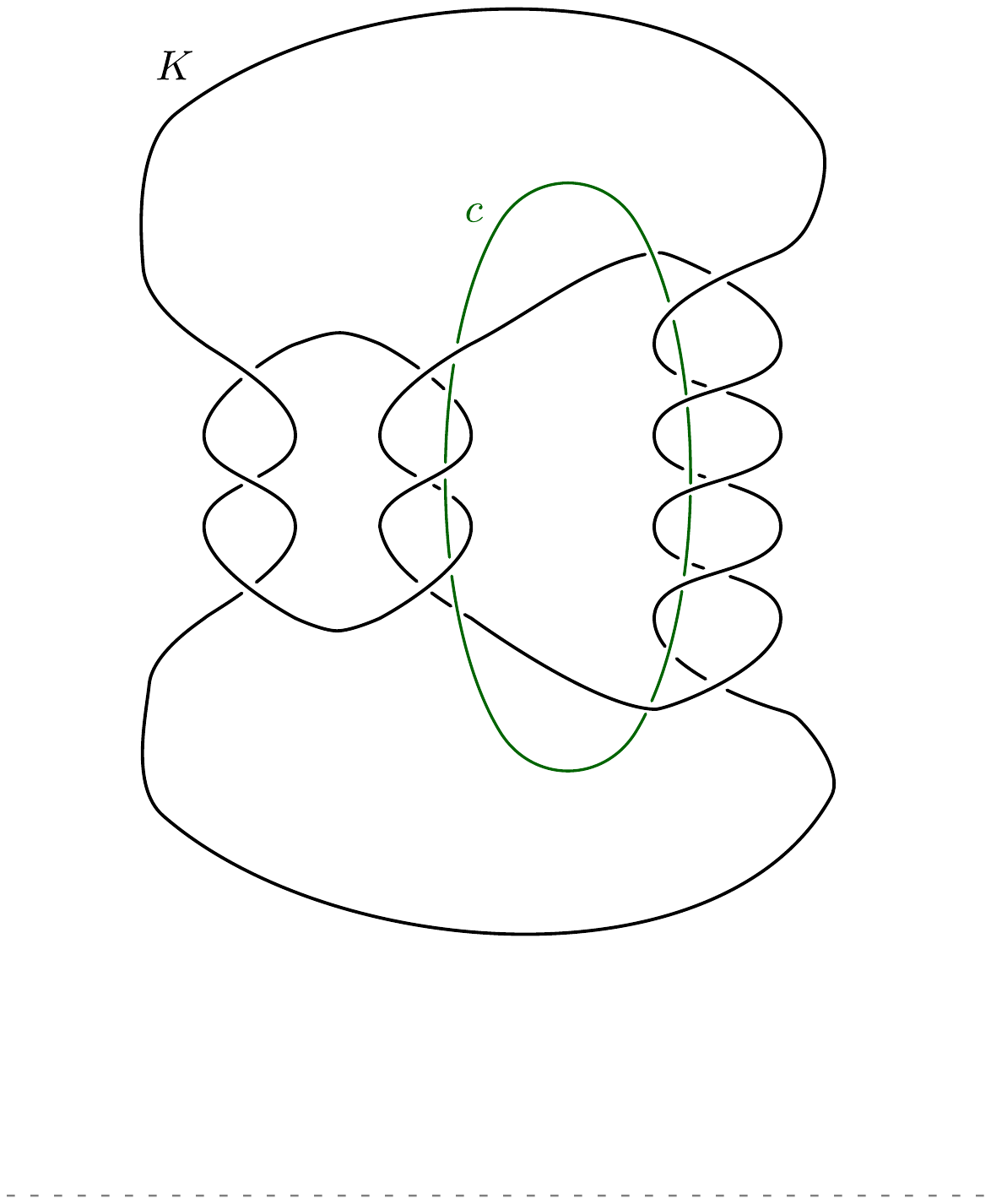}
\caption{Knot $K$ and a seiferter $c$ for it.}
\label{pretzel}
\end{figure}

Suppose we do a surgery on $c$ with the slope given by the ordinary fibres in the boundary of its fibred solid torus neighbourhood. We get a lens space in which we denote by $K'$ the image of $K$. Now doing the surgery on $K'$ with the slope induced by the original slope of $1$ has the same effect as first doing $1$-surgery on $K$ and then doing the reducible surgery on the exceptional fibre of index $4$. This gives $L(5,-2)\#L(3,-1) = L(5,3)\#L(3,2)$.

We now want to know in what lens space does $K'$ live. For this we need to understand what surgery we performed on $c$, i.e. what curves are the fibres around it. Suppose these fibres are $(p_1,q_1)$-curves in the solid torus neighbourhood of $c$. We know that $p_1 = 4$ and since there is a Seifert surface for $K$ that does not intersect $c$, the linking number of $K$ and $c$ is 0, hence by equation \eqref{excep} we find
$$
q_1 = \frac{0\cdot 4 + 15}{1} = 15.
$$

In the solid torus complementary to an open solid torus neighbourhood of $c$ the fibres on the boundary are thus $(15, 4)$-curves and therefore the surgery on $c$ gives the lens space $L(15, 4)$.

Putting this all together, we have a null-homologous genus $1$ knot $K'$ in $L(15,4)$ with an integral surgery giving $L(5,3)\#L(3,2)$. We now want to show that the knot $K' \subset L(15,4)$ is hyperbolic.

Since in \cite{bakerCabling} Baker has proven his conjecture for non-hyperbolic knots, we only need to consider the following cases:

\

\begin{itemize}
\item[\bfseries Case 1: \mdseries] $K'$ is contained in a ball;
\item[\bfseries Case 2: \mdseries] $K'$ is a Klein bottle knot and the slope is given by the surface slope of the essential annulus in its exterior;
\item[\bfseries Case 3: \mdseries] $K'$ is a torus knot and the slope is given by the surface slope of the essential annulus in its exterior;
\item[\bfseries Case 4: \mdseries] $K'$ is a cable knot and the slope is given by the surface slope of the essential annulus in its exterior.
\end{itemize}

\

\bfseries Case 1: \mdseries If a knot $K$ in $L(p,q)$ is contained in a ball, then all surgeries on it always have an $L(p,q)$ summand. This is clearly not true in our example.

\bfseries Case 2: \mdseries As the only lens spaces containing Klein bottle knots are of the form $L(4k,2k-1)$ \cite{bredonWood}, we see that this case is ruled out too.

\bfseries Case 3: \mdseries We can always fibre the lens space in such a way that the (non-trivial) torus knot we are considering is an ordinary fibre and the reducible slope is along the other fibres around it. Lens spaces can be represented as Seifert fibred spaces over $S^2$ with $0$, $1$ or $2$ exceptional fibres \cite{hatcher3mflds} but only the last case will result in a non-trivial connected sum of lens spaces. So suppose $L(15, 4) = S^2((p_1, x_1),(p_2,x_2))$. Then removing an ordinary fibre and filling it in the reducible way gives $L(p_1, x_1)\#L(p_2,x_2)$ (in the last expression $x$'s are only defined modulo the corresponding $p$'s). But also $S^2((p_1, x_1),(p_2,x_2))$ is the lens space with the order of the first homology $|p_1x_2+x_1p_2|$. In our case $p_1 = 3$, $p_2 = 5$, $x_1 = 2 + 3m$ and $x_2 = 3 + 5n$ for some integers $m$, $n$. So 
$$
\pm 15 = p_1x_2+x_1p_2 = 19 + 15(m+n),
$$

which clearly leads to a contradiction.

\bfseries Case 4: \mdseries Suppose $K'$ is a cable knot. Then its companion has a non-integral lens space surgery, so it has Seifert fibred exterior by the Cyclic Surgery Theorem \cite{CGLS} and hence is a torus knot by \cite[Theorem 6.1]{bakerBuck}. 

Let $W$ be a Heegaard solid torus of $L(15,4)$ and isotope $K'$ into $W$ in such a way that there is an identification of $W$ with a Heegaard solid torus of $S^3$ after which $K'$ becomes a cable of a torus knot in $S^3$. We fix the longitude of the companion torus knot using this identification.

Using this longitude, let $K'$ be a $(p,q)$-cable of an $(r,s)$-torus knot $T$ in \linebreak $W \subset L(15,4)$. Then the result of the reducible surgery on $K'$ is equal to the connected sum of $L(p,q)$ and the $q/p$-surgery on $T$. To get a $q/p$-surgery on $T$ we may first perform this surgery in $W$ and then attach a solid torus to the boundary of $W$. Let $L$ and $M$ be respectively the longitude and the meridian of $W$ (fixed by the same identification with a Heegaard solid torus of $S^3$). Then to obtain $L(15,4)$ we glue a solid torus $V$ to $W$ in such a way that the meridian of $V$ becomes a curve given by $15L+xM$ for some integer $x$.

In homology we have $l = rL$ and $M = rm$. Thus the first homology of the $q/p$-surgery on $T$ has a presentation matrix

$$
\begin{pmatrix} q & pr \\ xr & 15 \\  \end{pmatrix},
$$

so the order of the first homology is $\pm(15q - pxr^2)$.

All in all, the orders of the two lens space summands have to be $\pm p$ and $\pm(15q-pxr^2)$. On the other hand, they have to be $5$ and $3$.

If $p=\pm 3$ then we must have $15q-pxr^2 = \pm 5$. But also in this case $p | (15q-pxr^2)$, i.e. $3 \ | \ 5$ -- a contradiction. The case $p = \pm 5$ is completely analogous.
\end{proof}

There exists another seiferter for the $1$-surgery on $K$ that becomes the same exceptional fibre as $c$ does. This is the asymmetric seiferter $c_1'$ of \cite[Lemma 7.5]{networking1}. Suppose $K''$ is the image of $K$ in the lens space obtained by surgery on $c_1'$ with the slope given by the ordinary fibres around $c_1'$. After doing a calculation similar to the one we performed above it is easy to see that $K'' \subset L(19,4)$ is a primitive knot (i.e. generates $H_1$) that gives $L(5,3)\# L(3,2)$ by surgery.

The same elementary method as we used for $K'$ fails to show that $K''$ is hyperbolic. SnapPy \cite{SnapPy}, however, does suggest that $K''$ is hyperbolic.

\section{Some spaces that cannot be obtained by surgery with a seiferter}
\label{sec:spaces_seiferters}

Theorem \ref{main2} gives an obstruction on Seifert fibred surgery with a seiferter. Conjecturally \cite{miyazakiMotegiIII} all Seifert fibred surgeries have seiferters, so if this conjecture is true, Proposition \ref{prop:seiferter_restrictions} provides an obstruction to Seifert fibred surgery in general. Alternatively, it could be used to disprove it. Currently, however, only one family of knots with small Seifert Fibred surgeries is not known to have seiferters -- the one obtained by Teragaito in \cite{teragaito}. Application of the tests from the proposition does not provide any interesting information in this case.

Proposition \ref{obstructExamples}, which we restate below, gives an infinite family of small Seifert Fibred spaces that are obstructed from being surgeries with seiferters by Theorem \nolinebreak \ref{main2}.

{
\renewcommand{\thetheorem}{\ref{obstructExamples}}
\begin{prop}
Let $p \equiv 3 \ (\mbox{mod } 17)$. Then $S^2((p,-17),(2p-1,17),(2p+1,17))$ cannot be obtained by a Seifert surgery on a knot in $S^3$ with a seiferter.
\end{prop}
\addtocounter{theorem}{-1}
}

\begin{proof}
We calculate $H = 17$, thus we need to satisfy the obstructions of equation \eqref{obstructRatHom}. Since $p$, $2p-1$ and $2p+1$ are all coprime to $17$ and $-1$ is a quadratic residue modulo $17$, the obstructions of \eqref{obstructRatHom} are equivalent to none of
$$
(4p^2-1)p^*, (2p^2+p)(2p-1)^*, (2p^2-p)(2p+1)^* \mbox{ and } 4p^3-p
$$

being squares modulo $17$, where by $x^*$ we denote the multiplicative inverse of $x$ modulo $17$.

If $p \equiv 3 \ (\mbox{mod } 17)$ then
\begin{itemize}
\item $(4p^2-1)p^* \equiv 6 \ (\mbox{mod } 17)$;
\item $(2p^2+p)(2p-1)^* \equiv -6 \ (\mbox{mod } 17)$;
\item $(2p^2-p)(2p+1)^* \equiv 7 \ (\mbox{mod } 17)$;
\item $4p^3-p \equiv 3 \ (\mbox{mod } 17)$.
\end{itemize}

None of these are squares modulo $17$.
\end{proof}

We have also performed a simple computer search for spaces that are obstructed from being surgeries with seiferters by Theorem \ref{main2}. Out of $41468$ small Seifert fibred spaces with cyclic first homology that we checked, $17994$ spaces that satisfy the obstruction (i.e. cannot be obtained by a surgery with a seiferter) were found.

One such example with the lowest order of the first homology we were able to find (and one of the lowest one could hope for) is $S^2((2,-3),(3,1),(7,9))$ -- its first homology is $\mathbb{Z}_5$. However, due to Heegaard Floer homological reasons this particular space cannot be obtained by integral surgery on a knot in $S^3$. More concretely, the $d$-invariants of this space are $0, -2/5, -2/5, -8/5$ and $-8/5$. If it were obtained by integral (thus $\pm 5$) surgery on a knot in $S^3$, then its $d$-invariants (in some order) would differ by even integers from those of one of $L(5,\pm 1)$ (see \cite{NiWu}). However, one of their $d$-invariants is $\pm 1$. We remark that this space \slshape is \upshape a result of a non-integral surgery on a trefoil.

\section{Some observations about seiferters}
\label{sec:observations}

In Section \ref{sec:overview} we described some cases of the following situation. Given a hyperbolic knot $K$ with a Seifert fibred surgery with a seiferter $c$, there was a surgery on $c$, such that the image of $K$ in the resulting lens space was no longer hyperbolic. We were interested in a particular surgery on $c$, but it turns out that in very many cases this still holds, even though we might do a surgery on $c$ not with the slope given by preimages of the ordinary fibres around it.

In fact, we succeeded with finding such a surgery on $c$ in almost all cases we tried. In particular, this holds for: all Berge knots; all knots of Eudave-Mu\~noz from \cite{eudave-munoz}; first of the two examples in \cite{mattmanMiyazakiMotegi}; the non-symmetric seiferter $c'_1$ from \cite[Chapter 7]{networking1}. In all these examples the non-hyperbolic knots that we got were unknots, torus knots or cables of torus knots.

One technique of demonstrating such surgeries is via the Montesinos trick, as in the example we gave in Section \ref{sec:overview}. Sometimes one can also see such surgery directly from the link diagram. Recall that if we do a surgery on one component of a link we can handleslide the other components in the resulting space using the surgery slope. In other words, given a link $K \cup c$, if we perform a surgery on $c$, the knot $K'$ which is the image of $K$ under the surgery is isotopic to any other knot that can be obtained from $K$ by taking repeated band sums of $K$ with the surgery slope on $c$.

\begin{figure}
\includegraphics[scale=0.3, clip = true, trim = 90 160 110 150]{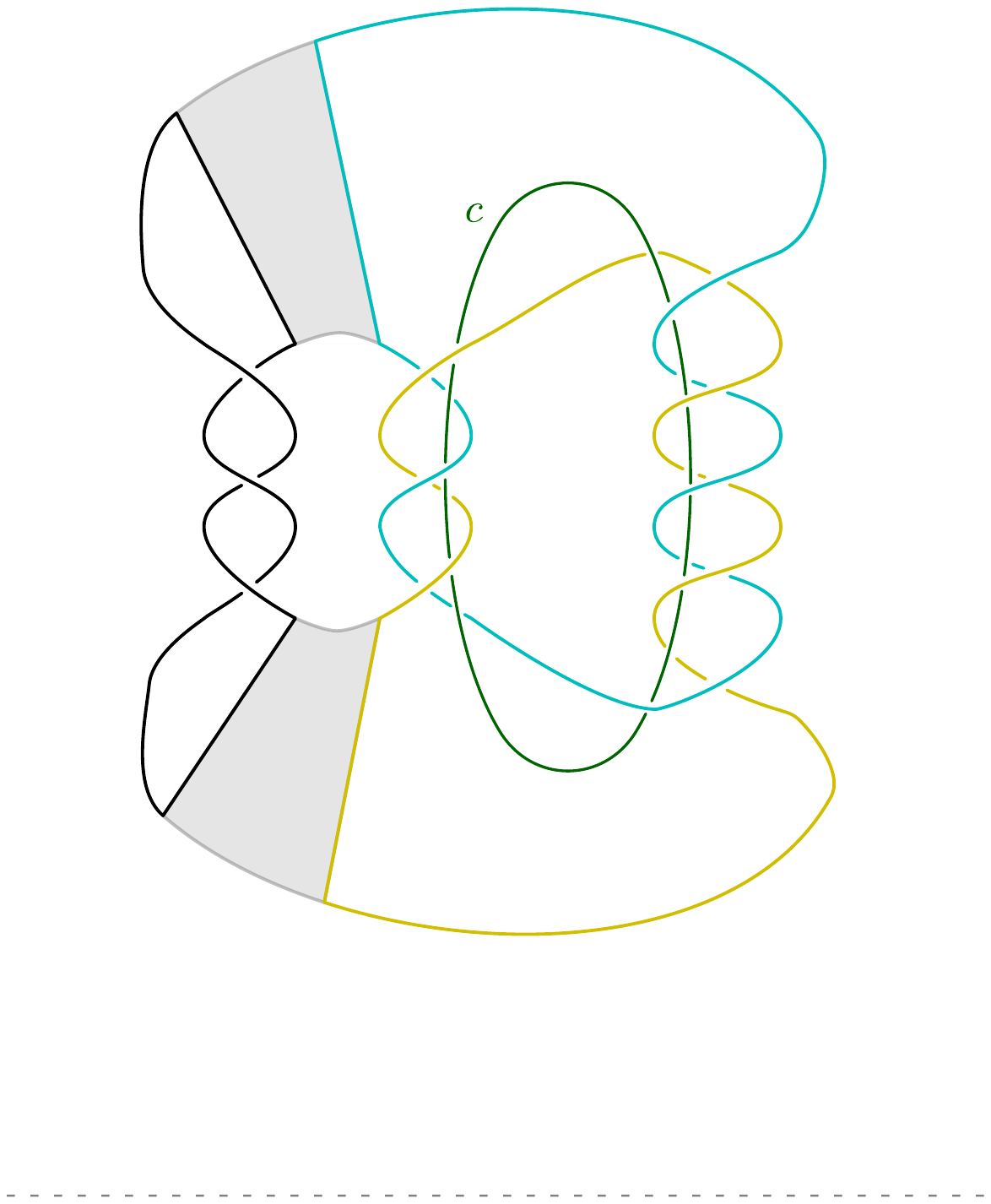}
\caption{The knot of Figure \ref{pretzel} can be decomposed as two band sums (bands are shown in grey) of the unknot (in black) with parallel cables of $c$ (in blue and yellow). This means that after some surgery on $c$ the knot becomes trivial in the resulting lens space.}
\label{pretzel2}
\end{figure}

\begin{figure}
\includegraphics[scale=0.3, clip = true, trim = 90 135 110 150]{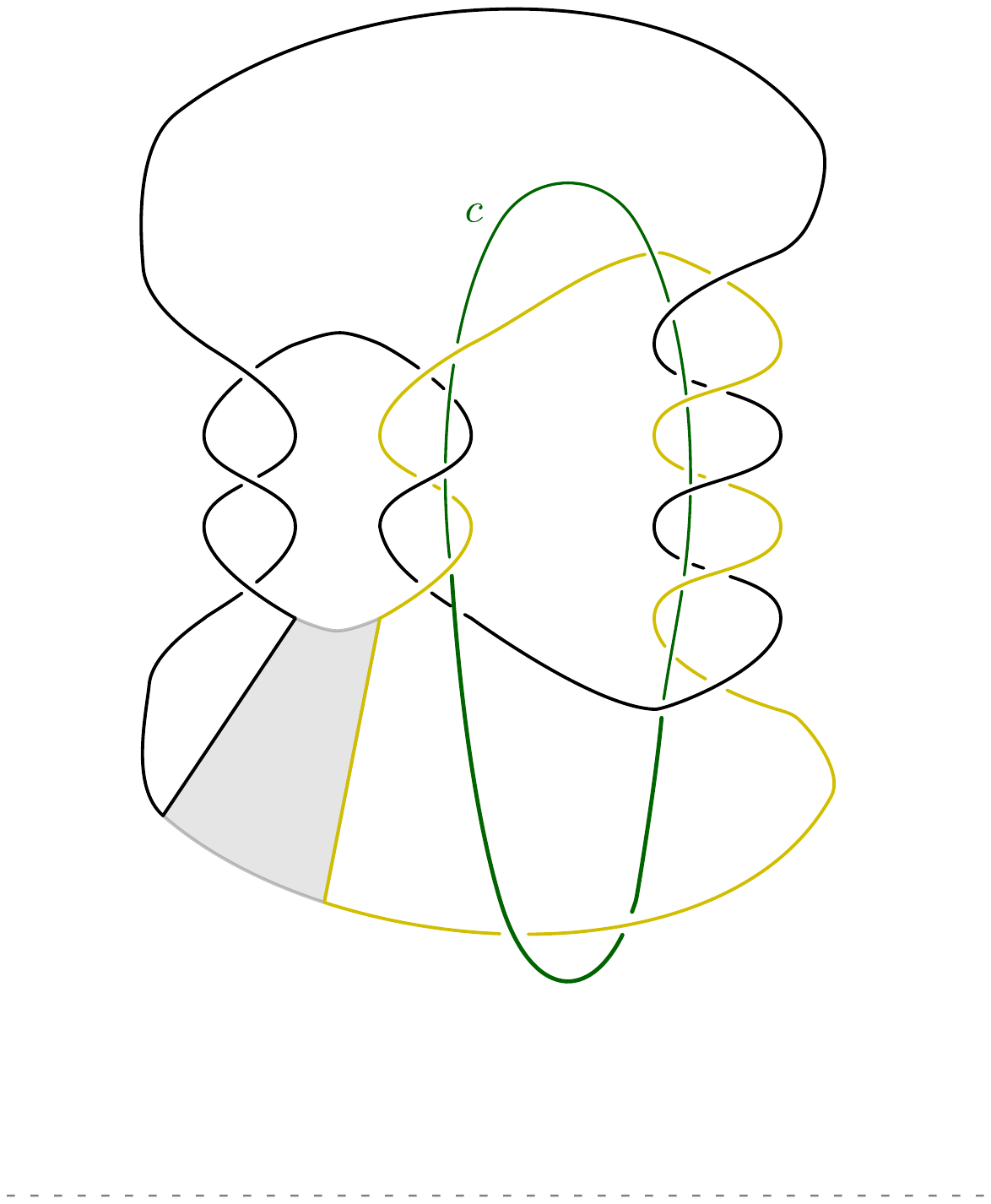}
\caption{Here $c$ is the asymmetric seiferter of \cite[Lemma 7.5]{networking1}. The knot becomes a torus knot after some surgery on $c$.}
\label{pretzel3}
\end{figure}

For example, the knot we used for our construction (see Figure \ref{pretzel}) can be decomposed as in Figure \ref{pretzel2}. This shows that there is a surgery on $c$ that transforms $K$ into the unknot. A similar argument works for the asymmetric seiferter for the same knot which becomes a torus knot in some lens space -- see Figure \ref{pretzel3}.

This makes us curious to know the answer to the following

\begin{qn}
Let $K$ be a hyperbolic knot with a Seifert fibred surgery with a seiferter $c$. Does there always exist a surgery on $c$ such that the image of $K$ in the resulting space is no longer hyperbolic?
\end{qn}

\begin{figure}
\includegraphics[scale=0.3, clip = true, trim = 85 160 105 150]{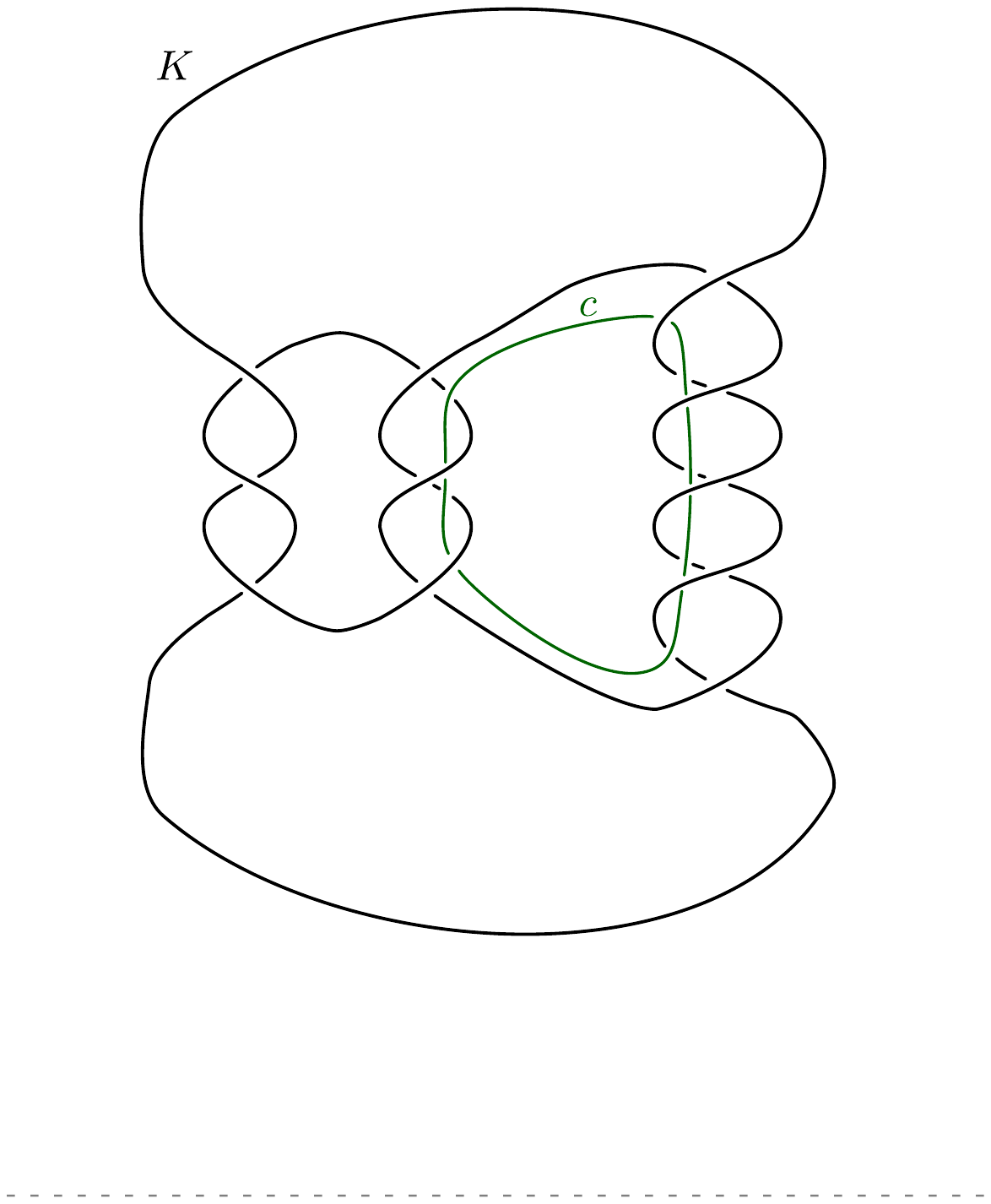}
\caption{$2$-surgery on $c$ makes $K$ a torus knot.}
\label{pretzel4}
\end{figure}

It is not immediately clear that there exists a filling of a seiferter on Figure \ref{pretzel4} that turns the knot into a non-hyperbolic one. However, using SnapPy Ken Baker observed that $2$-filling does have this property -- the fundamental group of the resulting knot exterior is $<a,b \ | \ a^2 = b^4>$, so the knot $K$ becomes a torus knot.

Observations about the Berge knots allow us to prove the following proposition. It seems somewhat similar to the results of \cite{networking2} but we do not know of any formal relation between these different viewpoints.

{
\renewcommand{\thetheorem}{\ref{berge}}
\begin{prop}
Every Berge knot can be obtained from a torus knot or a cable of a torus knot by repeatedly taking band sums with a cable of some unknot. For each Berge knot the unknot and its cable are fixed.
\end{prop}
\addtocounter{theorem}{-1}
}

\begin{proof}
Suppose we have a situation as in these examples, i.e. upon doing some surgery on a seiferter the knot becomes a torus knot or a cable of a torus knot. We can isotope it to lie in a standard way in one of the Heegaard solid tori. Isotoping in this lens space corresponds to the usual isotopies in the solid torus and band sums with the meridian of the other Heegaard solid torus \cite[Proposition 2.19(1)]{networking1}. Thus all examples we considered must be obtained by taking band sums with the same slope on a single unknot from torus knots or cables of torus knots.

Now we only need to demonstrate that for all Berge knots there is an unknot such that after a surgery on this unknot the Berge knot we consider becomes a torus knot or a cable of a torus knot. We have already done so for Berge knots of types IX-X.

Berge knots of types I-II are themselves torus knots and cables of torus knots. As shown in \cite[Proposition 7.2]{networking2} Berge knots of types VII-VIII have complexities at most $1$ (i.e. can be obtained from torus knots by twisting a number of times along a fixed unknot), so they also satisfy the property we are proving.

Therefore we only need to verify this for Berge knots of types III-VI and XI-XII. All of these have a surgery description on the minimally twisted five chain link and the tangles that give them via double branched covering are given in Figures 37-40 and 42 of \cite{bakerBerge2} (note that we use tangle conventions of \cite{networking4}).

Our proof by pictures will proceed as follows. We will reproduce the figures from \cite{bakerBerge2} with one additional tangle excised. Our tangles will thus have two boundary spheres, on projections denoted by a red and a green circle. The red circle will correspond to a Berge knot $K$ and the green circle $c$ to an unknot (a seiferter, but we don't need to know this) such that a surgery on it makes the original knot a torus knot (we choose this to include the core of one of the solid tori) or a cable of a torus knot.

The meridional surgery on $K$ will always correspond to the $\infty$-filling of the red tangle. Thus to verify that the green tangle indeed corresponds to an unknot one needs to see that after the $\infty$-filling of the red tangle the resulting tangle is homeomorphic to the trivial one (i.e. is a rational tangle). We will not depict this step, it is usually fairly straightforward to verify.

The surgery on $c$ that will make (the image of) $K$ a torus knot or a cable of a torus knot will always correspond to the $\infty$-filling of the green tangle. To verify that it indeed gives us a torus knot or a cable of a torus knot we need to demonstrate that the resulting tangle is homeomorphic to either a rational tangle, a sum of two rational tangles or a union of a rational tangle and a tangle corresponding to a cable space (as in Figures \ref{bergeSporadicAB3} and \ref{bergeSporadicAB4}). This is the step we provide pictures for.

\begin{figure}[h!]
\includegraphics[scale=0.5, clip = true, trim = 5 370 95 170]{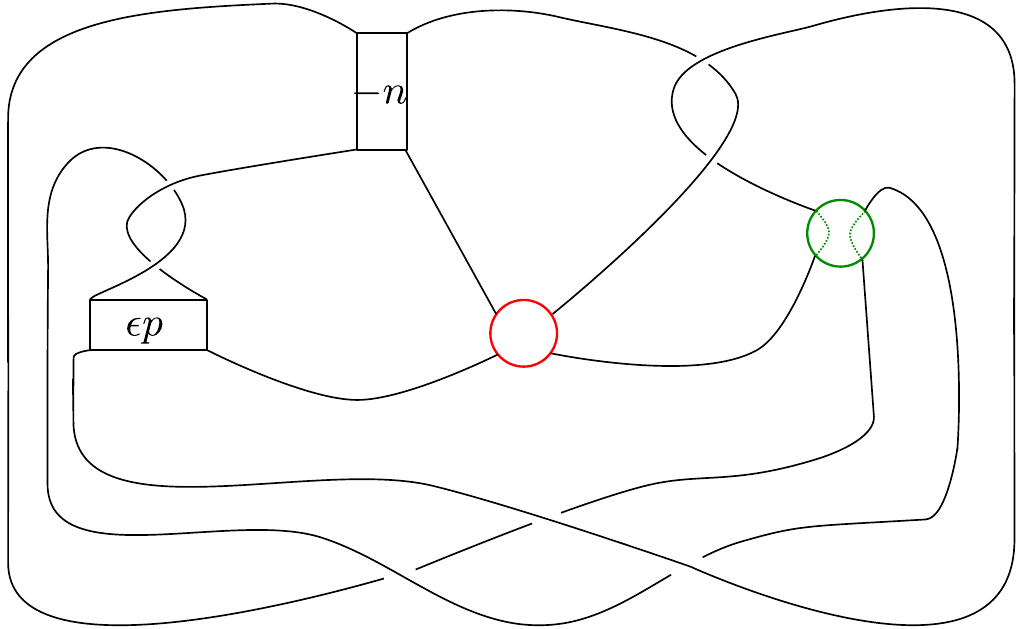}
\caption{Tangle description for the link formed by a Berge knot of type III (red) and an unknot (green). This corresponds to Figure 37 of \cite{bakerBerge2}.}
\label{bergeTypeIII1}
\end{figure}

\begin{figure}[h!]
\includegraphics[scale=0.5, clip = true, trim = 30 405 160 120]{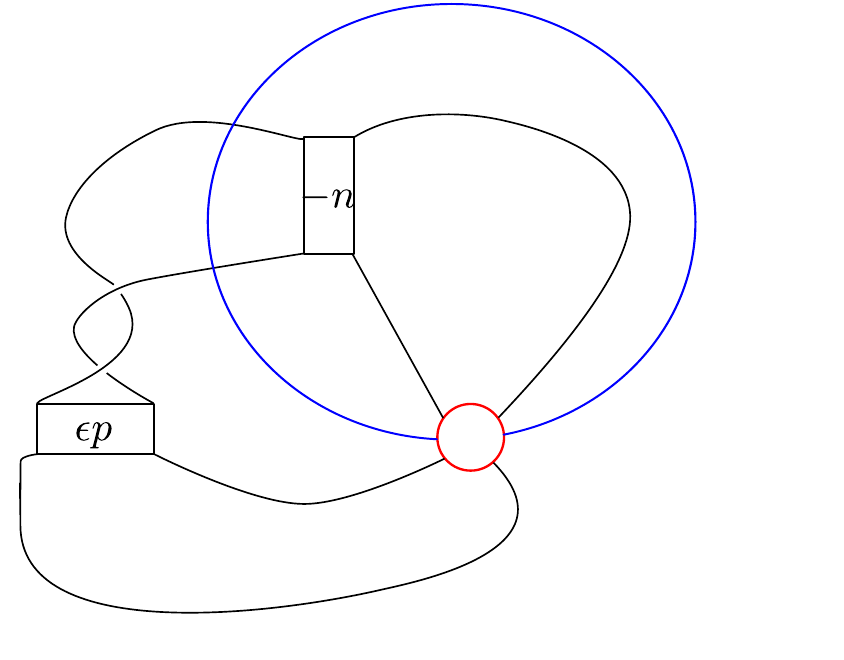}
\caption{Tangle description of the image of a Berge knot of type III after the $\infty$-filling of the green tangle from Figure \ref{bergeTypeIII1}. The knot becomes a torus knot, the blue disc lifts to the essential annulus that separates two solid tori.}
\label{bergeTypeIII2}
\end{figure}

\begin{figure}[h!]
\includegraphics[scale=0.5, clip = true, trim = 5 370 95 160]{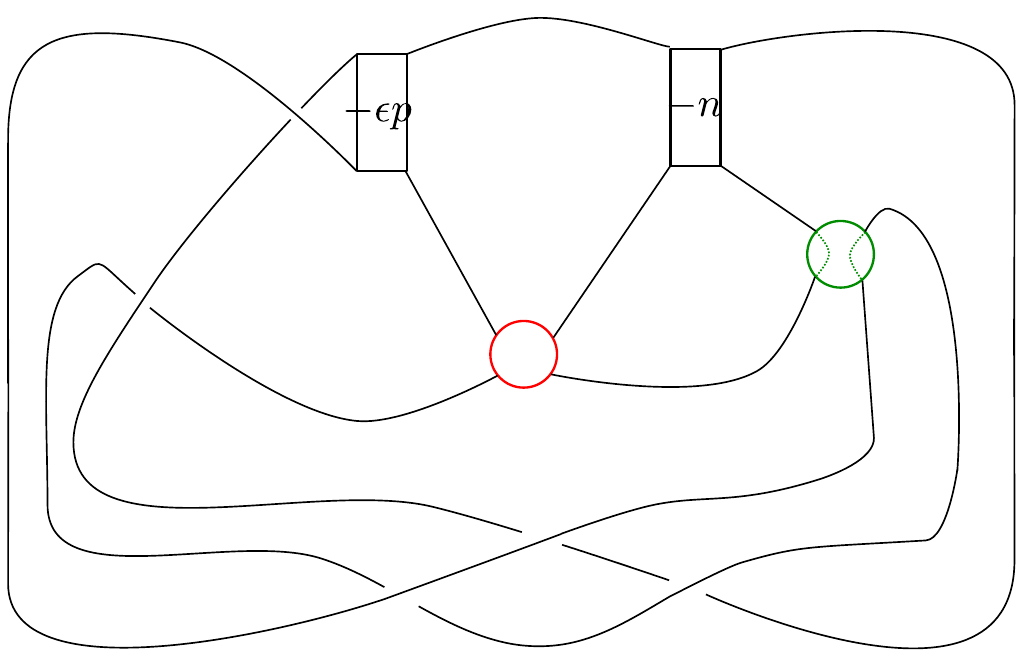}
\caption{Tangle description for the link formed by a Berge knot of type IV (red) and an unknot (green). This corresponds to Figure 38 of \cite{bakerBerge2}.}
\label{bergeTypeIV1}
\end{figure}

\begin{figure}[h!]
\includegraphics[scale=0.5, clip = true, trim = 0 370 80 160]{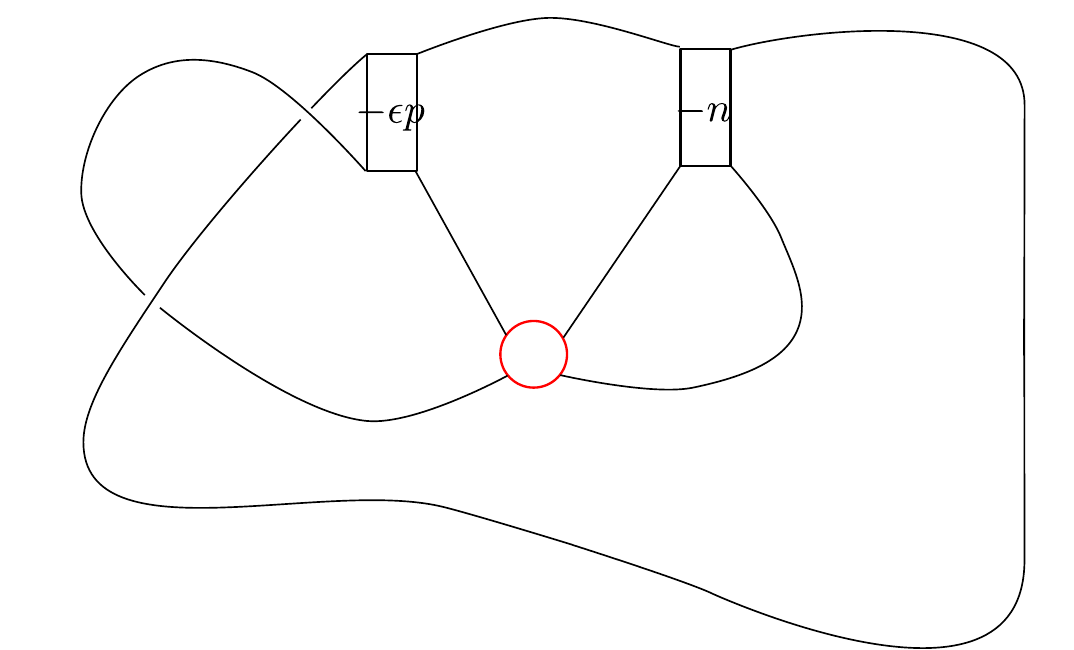}
\caption{Tangle description of the image of a Berge knot of type IV after the $\infty$-filling of the green tangle from Figure \ref{bergeTypeIV1}. The knot becomes a core of one of the Heegaard solid tori since the tangle is rational.}
\label{bergeTypeIV2}
\end{figure}

\begin{figure}[h!]
\includegraphics[scale=0.5, clip = true, trim = 5 370 95 170]{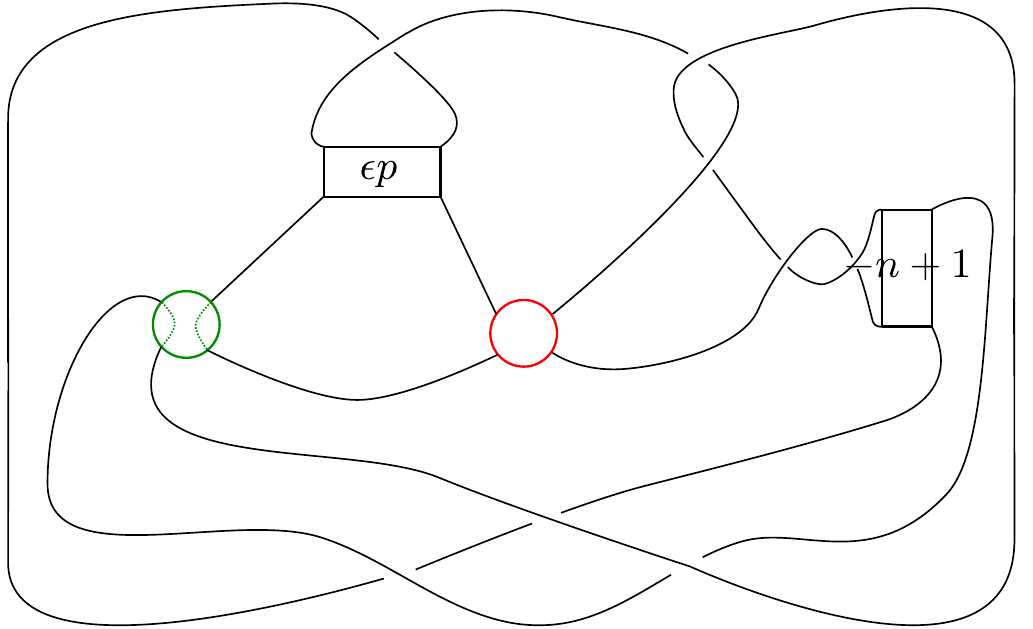}
\caption{Tangle description for the link formed by a Berge knot of type V (red) and an unknot (green). This corresponds to Figure 39 of \cite{bakerBerge2}.}
\label{bergeTypeV1}
\end{figure}

\begin{figure}[h!]
\includegraphics[scale=0.5, clip = true, trim = 0 370 40 120]{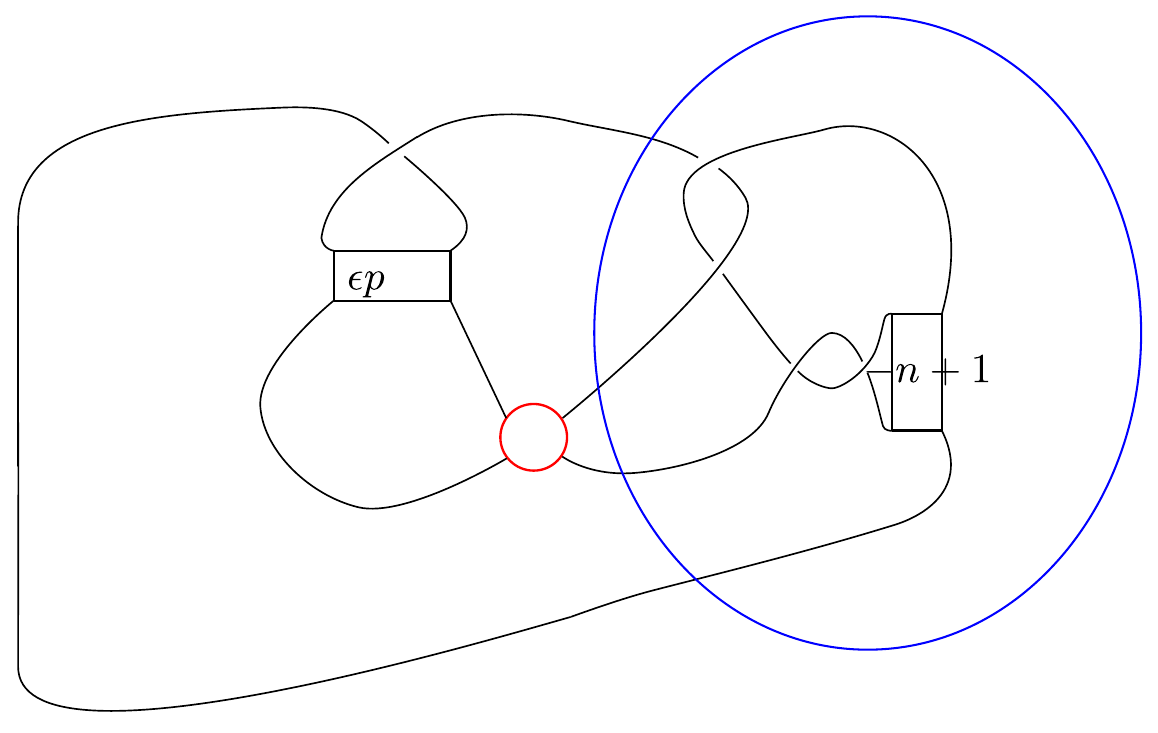}
\caption{Tangle description of the image of a Berge knot of type V after the $\infty$-filling of the green tangle from Figure \ref{bergeTypeV1}. The knot becomes a cable of a torus knot, the blue sphere lifts to the essential torus that separates the cable space from the torus knot exterior.}
\label{bergeTypeV2}
\end{figure}

\begin{figure}[h!]
\includegraphics[scale=0.5, clip = true, trim = 5 370 95 160]{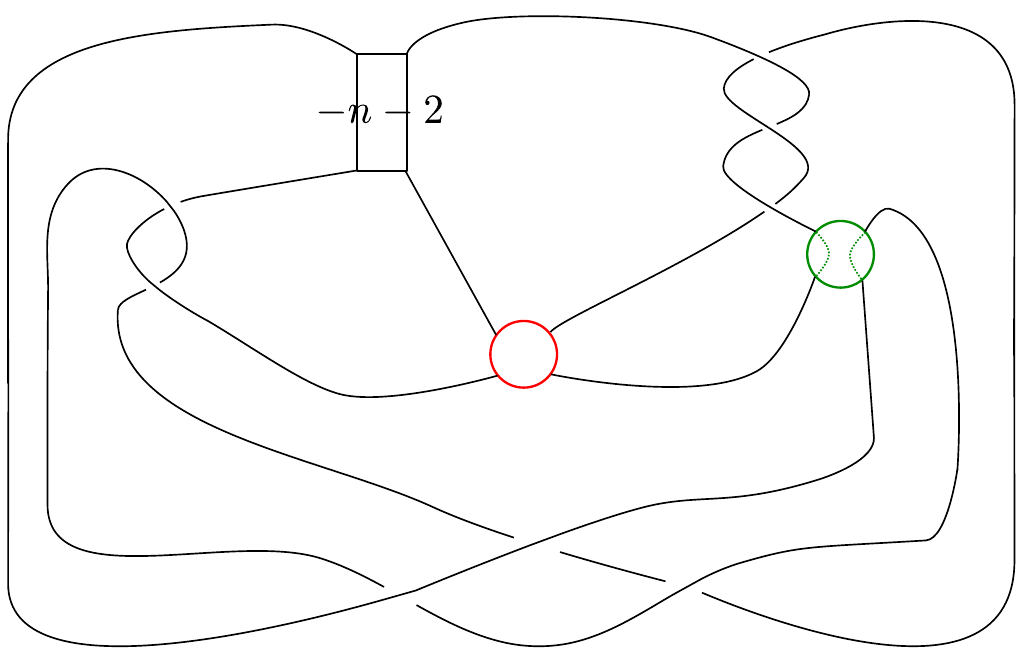}
\caption{Tangle description for the link formed by a Berge knot of type VI (red) and an unknot (green). This corresponds to Figure 40 of \cite{bakerBerge2}.}
\label{bergeTypeVI1}
\end{figure}

\begin{figure}[h!]
\includegraphics[scale=0.5, clip = true, trim = 0 370 80 160]{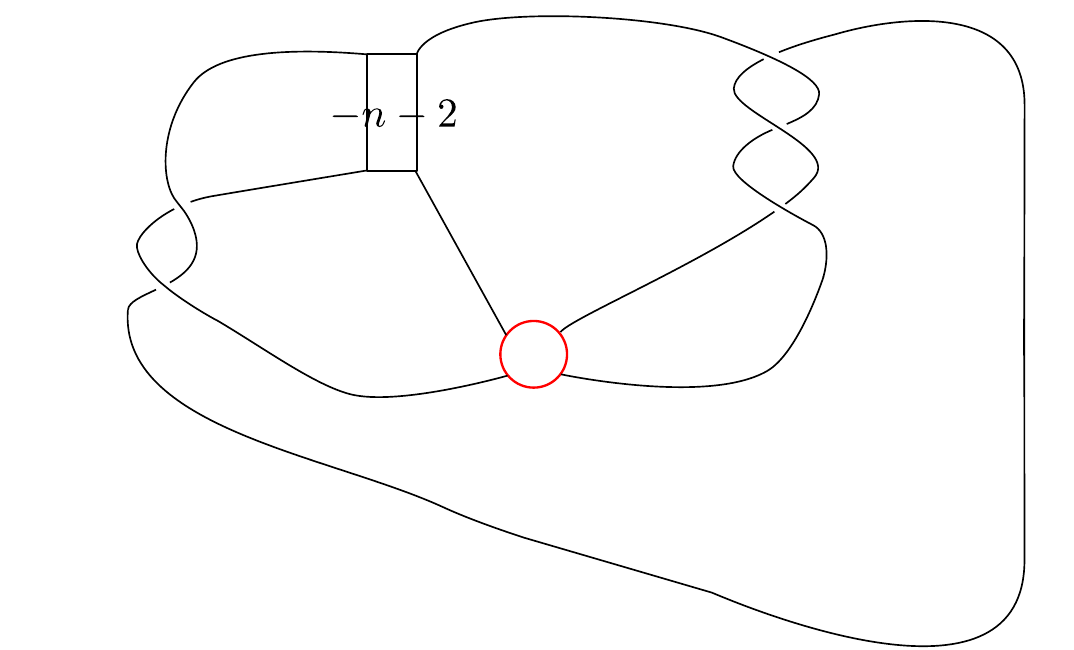}
\caption{Tangle description of the image of a Berge knot of type VI after the $\infty$-filling of the green tangle from Figure \ref{bergeTypeVI1}. The knot becomes a core of one of the Heegaard solid tori since the tangle is rational.}
\label{bergeTypeVI2}
\end{figure}

\begin{figure}[h!]
\includegraphics[scale=0.5, clip = true, trim = 5 370 95 170]{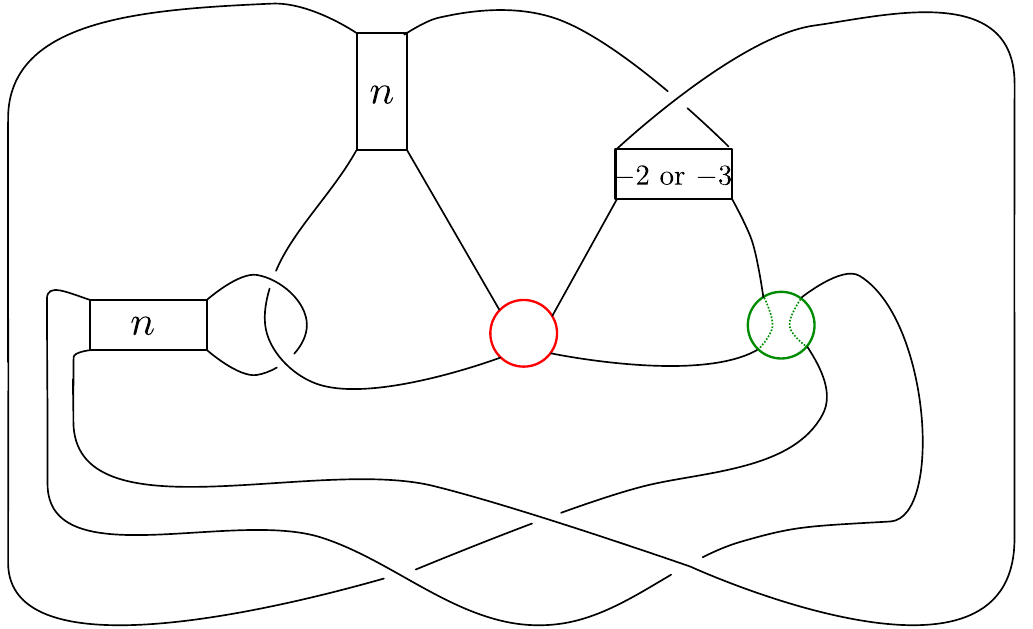}
\caption{Tangle description for the link formed by a Berge knot of type XI or XII (red) and an unknot (green). This corresponds to Figure 42 of \cite{bakerBerge2}.}
\label{bergeSporadicCD1}
\end{figure}

\begin{figure}[h!]
\includegraphics[scale=0.5, clip = true, trim = 30 395 140 160]{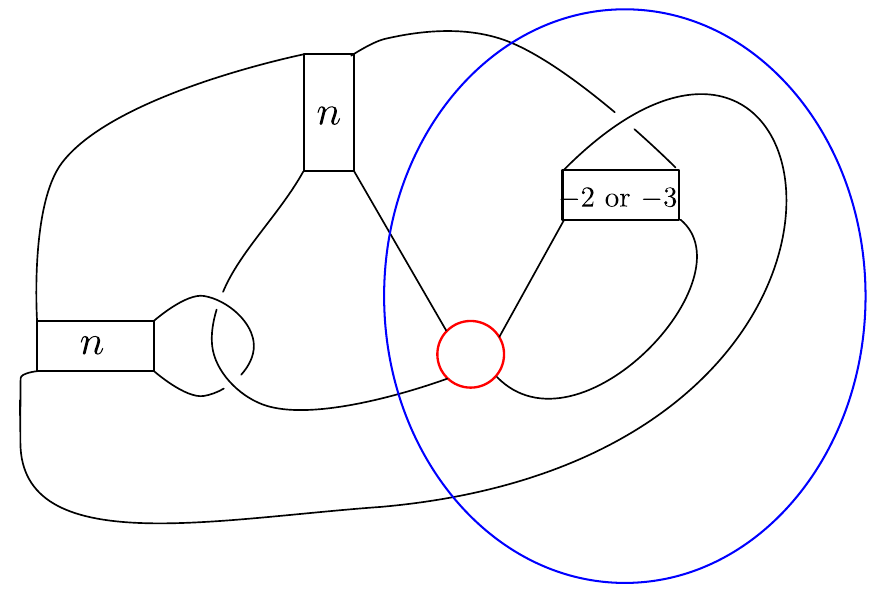}
\caption{Tangle description of the image of a Berge knot of type XI or XII after the $\infty$-filling of the green tangle from Figure \ref{bergeSporadicCD1}. The knot becomes a cable of a torus knot, the blue sphere lifts to the essential torus.}
\label{bergeSporadicCD3}
\end{figure}

This sequence of figures finishes the proof.
\end{proof}

\bibliographystyle{plain}  
\bibliography{biblio}

\end{document}